\documentclass[a4paper,12pt,reqno]{amsart}
\usepackage{amssymb}
\usepackage[dvips]{graphicx}

\input epsf.tex
\newdimen\xsize
\newdimen\oldbaselineskip
\newdimen\oldlineskiplimit
\def\restorelineskip{\baselineskip=\oldbaselineskip%
\lineskiplimit=\oldlineskiplimit}
\def\putm[#1][#2]#3{
\hbox{\vbox to 0pt{\parindent=0pt%
\vskip#2\xsize\hbox to0pt{\hskip#1\xsize $#3$\hss}\vss}}}%
\long\def\Line#1{\hbox to \hsize{#1}}
\def\putt[#1][#2]#3{
\vbox to 0pt{\noindent\hskip#1\xsize\lower#2\xsize%
\vtop{\restorelineskip#3}\vss}}

\makeatletter
\def\xbig[#1]#2{{\hbox{$\m@th\left#2\vbox to#1\xsize{}%
\right.\n@space$}}}
\makeatother
\def\xlar[#1]#2{%
\smash{\mathop{ \hbox to #1\xsize{\leftarrowfill}}\limits^{#2}}}
\def\xrar[#1]#2{%
\smash{\mathop{ \hbox to #1\xsize{\rightarrowfill}}\limits^{#2}}}
\def\xline[#1]{\hbox to #1\xsize{\leaders\hrule\hfill}}

\DeclareFontFamily{U}{rsf}{\skewchar\font'177}%
\DeclareFontShape{U}{rsf}{m}{n}{<-6>rsfs5<6-8>rsfs7<8->rsfs10}{}%
\DeclareFontShape{U}{rsf}{b}{n}{<-6>rsfs5<6-8>rsfs7<8->rsfs10}{}%
\DeclareMathAlphabet\RSFS{U}{rsf}{m}{n}
\SetMathAlphabet\RSFS{bold}{U}{rsf}{b}{n}
\def\sf#1{{\mathsf{#1}}}

\def\slsf{\slshape \sffamily }

\def\msmall#1{\mathchoice{\hbox{\small$\displaystyle {#1}$}}{#1}{#1}{#1}}

\def\cc{{\mathbb C}}

\def\pp{{\mathbb P}}

\def\lim{\mathop{\sf{lim}}}

\def\vol{\sf{Vol}}

\def\<{\langle}\let\la=\<
\def\>{\rangle}\let\ra=\>
  
\def\comp{\Subset}

\def\ddef{\mathrel{{=}\raise0.3pt\hbox{:}}}
\def\deff{\mathrel{\raise0.3pt\hbox{\rm:}{=}}}

\def\fraction#1/#2{\mathchoice{{\msmall{ #1\over#2}}}%
{{ #1\over #2 }}{{#1/#2}}{{#1/#2}}}

\def\longpoints{\leaders\hbox to 0.5em{\hss.\hss}\hfill \hskip0pt}
\def\stateskip{\smallskip}
\def\state#1. {\stateskip\noindent{\bf#1. }} 
\def\statep#1. {\stateskip\noindent{\bf#1 }} 
\def\proof{\state Proof. \2}

\def\Chi{\raise 2pt\hbox{$\chi$}}

\def\ie{\hskip1pt plus1pt{\sl i.e.\/,\ \hskip1pt plus1pt}}

\def\sli{{\sl i)} } 
\def\slii{{\sl i$\!$i)} }

\def\Chi{\raise 2pt\hbox{$\chi$}}

\let\phI=\phi\let\phi=\varphi\let\varphi=\phI
\let\cal=\mathcal

\def\calh{{\cal H}}

\def\calm{{\cal M}}

\def\comp{\Subset}

\def\1{{1\mkern-5mu{\rom l}}}

\def\ge{\geqslant}

\def\fraction#1/#2{\mathchoice{{\msmall{ #1\over#2}}}%
{{ #1\over #2 }}{{#1/#2}}{{#1/#2}}}

\newcommand{\2}{\thinspace}

\def\qed{\ \ \hfill\hbox to .1pt{}\hfill\hbox to .1pt{}\hfill $\square$\par}

\def\comment#1\endcomment{}

 \belowdisplayskip=\abovedisplayskip

\def\lineeqqno(#1){\hfill\llap{\vbox to 10pt%
{\vss\begin{align} \eqqno(#1)\end{align}\vss}}\vskip1pt}

\advance\headheight 1.2pt

\def\ShowwLLabel#1{}

\def\thechpt{\Roman{chpt}}

\def\newchapt[#1]#2{\newpage%
\refstepcounter{chpt}\setcounter{subsection}{0}%
\setcounter{thm}{0}\setcounter{defi}{0}%
\setcounter{rema}{0}\setcounter{exrc}{0}%
\renewcommand{\thesubsection}{\thechpt.\arabic{subsection}}%
\section*{\begin{center}\huge \bf Chapter \thechpt\\
#2 \end{center}}\label{#1}%
\ \smallskip%
\markboth{Chapter \thechpt}{#2}%
}
\def\newsect[#1]#2{\refstepcounter{section}\setcounter{equation}{0}%
\renewcommand{\thesubsection}{\arabic{section}.\arabic{subsection}}%
\section*{\arabic{section}.
#2}\vspace{-20pt}\label{#1}\vspace{20pt}%
\markboth{Section \arabic{section}}{#2}}

\def\newlect[#1]#2{\refstepcounter{section}%
\renewcommand{\thesubsection}{\arabic{section}.\arabic{subsection}}%
\section*{Lecture \arabic{section}\\
#2}\label{#1}%
\markboth{Lecture \arabic{section}}{#2}}

\def\newprg[#1]#2{\refstepcounter{subsection}%
\subsection*{{\thesubsection.\ #2}} \label{#1}%
}

\def\newappx[#1]#2{%
\refstepcounter{appx}\setcounter{section}{0}%
\renewcommand{\thesubsection}{A\arabic{appx}.\arabic{subsection}}%
\section*{Appendix \arabic{appx}\\ #2}
\label{#1}%
\markboth{Appendix A\arabic{appx}}{#2}
}

\newtheorem{thm}{Theorem}[section]
   \def\newthm#1{\begin{thm}\label{#1}}

\newtheorem{nnthm}{Theorem.} 
   \def\newnnthm#1{\begin{nnthm} \label{#1}}
\newtheorem{lem}{Lemma}[section]
   \def\newlemma#1{\begin{lem} \label{#1}}

\newtheorem{prop}{Proposition}[section]
   \def\newprop#1{\begin{prop}\label{#1}}

\newtheorem{corol}{Corollary}[section]
   \def\newcorol#1{\begin{corol} \label{#1}}

\newtheorem{defi}{Definition}[section]
   \def\newdefi#1{\begin{defi} \label{#1}\rm }

\newtheorem{exmp}{Example}[section]
   \def\newexmp#1{\begin{exmp} \label{#1}\rm }

\newtheorem{exrc}{Exercise}
   \def\newexrc#1{\begin{exrc} \label{#1}\rm }

\newtheorem{rema}{Remark}[section]
   \def\newrema#1{\begin{rema} \label{#1}\rm }

\def\eqqno(#1){\label{(#1)}}
\def\eqqref(#1){(\ref{(#1)})}

\title{On fixed points of rational self-maps of complex projective plane}
\author{S. Ivashkovich}
\date{\today}

\address{
Universit\'e de Lille-1, UFR de Math\'ematiques, 59655 Villeneuve
d'Ascq, France} \email{ivachkov@math.univ-lille1.fr}
\address{IAPMM Nat. Acad. Sci. Ukraine
Lviv, Naukova 3b, 79601 Ukraine}


\begin{document}
\begin{abstract}
For any given natural $d\ge 1$ we provide examples of rational self-maps of
complex projective plane  $\pp^2$ of degree $d$ without
(holomorphic) fixed points. This makes a contrast with the situation
in one dimension. We also prove that the set of fixed point free rational
self-maps of $\pp^2$ is closed (modulo ``degenerate'' maps) in some natural
topology on the space of rational self-maps of $\pp^2$.
\end{abstract}

\maketitle


\setcounter{tocdepth}{1}

\tableofcontents

\newsect[1]{Introduction} Along this note $\pp^2$ stands for the complex
projective plane as well as $\pp^1$ for the projective line. A
meromorphic self-map $f:\pp^2\to\pp^2$ can be viewed as a
holomorphic map $f:\pp^2\setminus \{ \text{ isolated points } \}\to
\pp^2$. Then the graph of $f$ extends to an analytic subvariety of
the product $\pp^2\times\pp^2$, this extension will be denoted as
$\Gamma_f$ and called the graph of the meromorphic mapping $f$. The
graph $\Gamma_f$ is then an algebraic subvariety by Chow's theorem
and consequently $f$ itself can be defined by a pair of rational
functions, \ie $f$ is rational.

\smallskip A {\slsf meromorphic} fixed point of $f$ is a point
$p\in\pp^2$ such that $p\in f[p]$. Here by $f[p]$ one means the full
image of $p$ by $f$:
\[
f[p]\deff \Gamma_f\cap \Big(\{p\}\times \pp^2\Big),
\]
\ie $f[p]$ can be a complex curve (such $p$ is, by definition, a
point of indeterminacy of $f$). By obvious homological reasons
$\Gamma_f$ intersects the diagonal $D$ in $\pp^2\times\pp^2$.
Therefore meromorphic fixed points for any $f:\pp^2\to\pp^2$ always
exist.

\smallskip In this paper we shall consider only {\slsf holomorphic}
fixed points: a point $p\in\pp^2$ is said to be a {\slsf
holomorphic} fixed point of $f$,  if $f$ is holomorphic in a
neighborhood of $p$ and $f(p)=p$. In what follows holomorphic fixed
points shall be called simply - fixed points. Let $X$ be a compact
complex manifold and $f:X\to X$ a meromorphic self-map. The
topological degree of $f$ is the number of preimages of a generic
point. The goal of this note is to prove the following:

\begin{thm}
\label{fixed} For any given integer $d\ge 1$ there exist rational
self-maps $f:\pp^2\to\pp^2$ of degree $d$ without holomorphic fixed
points.
\end{thm}

\smallskip One of the reasons for the interest in fixed points of
meromorphic maps lies in the attempt to understand what should be an
analog of a Lefschetz Fixed Point Formula in meromorphic case, see
ex. \cite{B} and \cite{BK}.

\medskip In Section 4 we define a natural topology on the space
$\calm (\pp^2)$ of meromorphic self-maps of the complex projective
plane. Denote by $FFix (\pp^2)$ the set of fixed point free
meromorphic self-maps of $\pp^2$. By $DFix(X)$ the subset of $\calm
(\pp^2)$ which consists of the maps with a {\slsf curve} of fixed
points (\ie in some sense they are degenerate maps). $DFix(X)$ is a
closed subset of $\calm (\pp^2)$. We prove the following:

\begin{thm}
\label{closed} The set $FFix(\pp^2)\cup DFix(\pp^2)$ is closed in
$\calm (\pp^2)$.
\end{thm}

\noindent I.e., a sequence of fixed point free meromorphic mapping
converge either to a fixed point free meromorphic map, or to a map
with a curve of fixed points. Notice that by the Theorem \ref{fixed}
the set $FFix(\pp^2)\cup DFix(\pp^2)$ is a proper subset of $\calm
(\pp^2)$.

\medskip\noindent{\slsf Acknowledgement.} 1. I heard of the question
whether any rational self-map of $\pp^2$ has a holomorphic fixed
point  for the first time in the talk of J.-E. Fornaess on the
``Colloque en l'honneur de P. Dolbeault'' in Paris, June 1992.

\smallskip\noindent 2. I would like also to give my thanks to the
organizers and participants of the Conference ``Dynamique et
G\'eom\'etrie Complexe'' in Luminy, June 2009, where this note was
basically done, for the encouraging atmosphere and stimulating
discussions. I am especially grateful to E. Bedford for sending me
his Lecture Notes \cite{B}.

\newsect[sect2]{Examples in the product of projective lines}

We shall work first with $\pp^1\times\pp^1$ as a model of the
complex projective plane, because in this manifold the geometric
picture is particularly clear and formulas are simpler. Then we
shall transfer our examples to $\pp^2$.

\newprg[prg2.1]{Equation for the fixed points}

From what was said above in the Introduction it is clear that in
order to produce mappings without fixed points one needs to produce
mappings $f$ such that $\Gamma_f$ intersects the diagonal $D$ in
$\big(\pp^1\times\pp^1\big)\times \big(\pp^1\times\pp^1\big)$ only
over the points of indeterminacy of $f$.

\smallskip
 Every rational self-map $f:\pp^1\times\pp^1\to
\pp^1\times\pp^1$ can be written as $f=(f_1,f_2)$, where
$f_j(z_1,z_2)$, $j=1,2$, are rational functions of two complex
variables. A point $p=(z_1,z_2)\in\pp^1\times\pp^1$ is a fixed point
for $f$ if and only if both $f_1$ and $f_2$ are holomorphic in a
neighborhood of $p$ and
\begin{equation}
\eqqno(fix-pnt)
\begin{cases}
f_1(z_1,z_2) = z_1,\cr
f_2(z_1,z_2) = z_2.
\end{cases}
\end{equation}
When writing \eqqref(fix-pnt) we mean always affine coordinates. If
a fixed point $p$ has one or both of its coordinates equal to
$\infty$ then one should take an appropriate affine coordinates in a
neighborhood of $p$ and appropriately rewrite $\eqqref(fix-pnt)$.
But holomorphicity here means holomorphicity {\slsf with values in
$\pp^1$}.

\begin{figure}[h]
\centering
\includegraphics[width=2.5in]{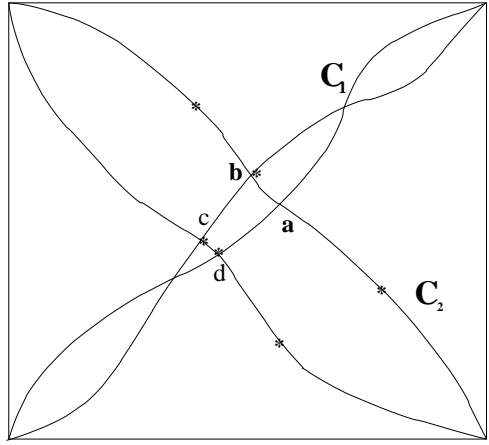}
\caption{Curves $C_1$ and $C_2$ are defined by equation \eqqref(fix-pnt).
Points of intersection of these curves are the fixed points of our map unless 
they are indeterminacy points at the same time. Points, marked by asterisks, 
are indeterminacy points. Therefore, on this picture $a$ is a fixed point, 
but $b$, $c$ and $d$ are not.}
\label{fix-fig}
\end{figure}

\smallskip
Set $C_1\deff\{ (z_1,z_2)\in\pp^1\times\pp^1: f_1(z_1,z_2)=z_1\}$
and $C_2\deff\{ (z_1,z_2)\in\pp^1\times\pp^1: f_2(z_1,z_2)=z_2\}$.
These are complex curves and their intersection $C_1\cap C_2$
contains all fixed points of $f$. More precisely, if by $I(f)$ we
denote the indeterminacy set of $f$ then $Fix(f) = (C_1\cap
C_2)\setminus I(f)$. Here $Fix(f)$ stands for the set of
(holomorphic) fixed points of $f$. Remark, finally, that a point
$p=(z_1,z_2)$ is an indeterminacy point of $f$ if it is an
indeterminacy point of at least one of $f_1$ or $f_2$.

\medskip

Let us try to find rational functions $f_2$ such that for
$f_1(z_1,z_2) = \frac{z_2}{z_1}$ the rational map $f=(f_1,f_2)\in
Rat(\pp^1\times\pp^1)$ has no fixed points. Write $f_2(z_1,z_2) =
\frac{P(z_1,z_2)}{Q(z_1,z_2)}$. Then $C_1=\{z_2=z_1^2\}$ and
$C_2=\{P(z_1,z_2)=z_2Q(z_1,z_2)\}$. Therefore
\begin{equation}
\eqqno(z-z2)
Fix(f) = \{ (z,z^2)\in\pp^1\times\pp^1: P(z,z^2)=z^2Q(z,z^2)\}
\setminus I(f).
\end{equation}
Remark that in order for $Fix(f)$ to be an empty set any solution
$\lambda$ of \eqqref(z-z2) should be either $0$ or $\infty$ or,
otherwise, $\lambda$ should be the root of both polynomials
$P(z,z^2)$ and $Q(z,z^2)$. More precisely the following is true:

\begin{lem}
\label{lem-z-z2} Let $P$ and $Q$ be relatively prime and suppose
that every non-zero root of $P(z,z^2)-z^2Q(z,z^2)$ is the root of
both $P(z,z^2)$ and $Q(z,z^2)$. Then the map
\[
f : (z_1,z_2) \to \Big( \frac{z_2}{z_1},\frac{P(z_1,z_2)}{Q(z_1,z_2)}\Big)
\]
has no fixed points in $\pp^1\times\pp^1$.
\end{lem}
\proof Suppose that $p=(z_1,z_2)$ is a fixed point of $f$. Then
$p\not= (0,0),(\infty,\infty)$ and $z_2=z_1^2$, \ie $p=(z,z^2)$ for
some non-zero complex number $z$. In addition we have that
$P(z_1,z_2)=z_2Q(z_1,z_2)$ and this implies that
$P(z,z^2)-z^2Q(z,z^2)=0$. Therefore $z$ is the root of both
$P(z,z^2)$ and $Q(z,z^2)$ by the assumption of the Lemma. But that
means that our point $p=(z,z^2)$ belong to the zero divisor of both
$P$ and $Q$. Since they are relatively prime, the point $p$ is an
indeterminacy point of $f_2=\frac{P}{Q}$.

\smallskip\qed

\begin{rema} \rm
Two polynomials $P(z_1,z_2)$ and $Q(z_1,z_2)$ are relatively prime
if their zero divisors do not have common component. This is easily
checked when $P$ and $Q$ are simple enough.
\end{rema}

\newprg[prg2.2]{Examples}

Now let us  give a list of examples following from Lemma \ref{lem-z-z2}.
Let's start, for the sake of clarity, with low degrees.

\begin{exmp} \rm
\label{exmp1}
Consider the map
\begin{equation}
f(z_1,z_2) = \Big( \frac{z_2}{z_1} , \frac{z_1^2-1}{z_2-1}\Big).
\end{equation}
\end{exmp}
In this case $P(z_1,z_2) = z_1^2-1$ and $Q(z_1,z_2) = z_2-1$. $P$
and $Q$ are obviously relatively prime. Moreover, polynomials
$P(z,z^2)= z^2-1$ and $Q(z,z^2)= z^2-1$ do satisfy the condition of
Lemma \ref{lem-z-z2}. Really: $z^2-1-z^2(z^2-1)= -(z^2-1)^2$.
Therefore this map has no fixed points. The degree of $f$ is $2$,
the indeterminacy points are: $(0,0), (\infty , \infty), (\pm 1,1)$.

\begin{exmp} \rm
\label{exmp2}
Consider the map
\begin{equation}
f(z_1,z_2) = \Big( \frac{z_2}{z_1} , \frac{(z_1-1)(z_1-2)(3z_1-2)}{z_2-3z_1 +2}\Big).
\end{equation}
\end{exmp}
In this case $P(z_1,z_2) = (z_1-1)(z_1-2)(3z_1-2)$ and
$Q(z_1,z_2) = z_2 - 3z_1 +2$. Therefore $P(z,z^2) = (z-1)(z-2)(3z-2)$,
$Q(z,z^2) = z^2 - 3z +2 = (z-1)(z-2)$ and $P(z,z^2)-z^2Q(z,z^2) =
(z-1)(z-2)(3z-2-z^2)=-(z-1)^2(z-2)^2$. Condition of
Lemma \ref{lem-z-z2} is again satisfied. The degree of $f$ is $3$, the
indeterminacy points are: $(0,0), (\infty , \infty),
(1,1), (2,4)$ and  $(\frac{2}{3}, 0)$.

\medskip Now let us give examples in all degrees.

\begin{exmp} {\slsf (Even degrees).} \rm
\label{exmp3}
Let $(a_1,b_1),...,(a_d,b_d)$ be such pairs of complex numbers that quadratic
polynomials $(z-a_1)(z-b_1),...,(z-a_d)(z-b_d)$
are pairvise relatively prime. Consider the rational map
\begin{equation}
f(z_1,z_2) = \Big( \frac{z_2}{z_1} , \frac{(z_1^2-1)\prod_{j=1}^d(z_1-a_j)
(z_1-b_j)}{\prod_{j=1}^d[z_2-(a_j+b_j)z_1 + a_jb_j]}\Big).
\end{equation}
Polynomials $P(z_1,z_2) = (z_1^2-1)\prod_{j=1}^d(z_1-a_j)(z_1-b_j)$ and
$Q(z_1,z_2) = \prod_{j=1}^d[z_2-(a_j+b_j)z_1 + a_jb_j]$
are obviously relatively prime (think about their zero divisors). To
check Lemma \ref{lem-z-z2} write
\[
P(z,z^2) - z^2Q(z,z^2) = (z^2-1)\prod_{j=1}^d(z-a_j)(z-b_j) -
z^2\prod_{j=1}^d[z^2-(a_j+b_j)z + a_jb_j] =
\]
\[
= \prod_{j=1}^d(z-a_j)(z-b_j)(z^2-1-z^2) = - \prod_{j=1}^d(z-a_j)(z-b_j).
\]
Therefore the condition of Lemma \ref{lem-z-z2} is satisfied. The degree
of $f$ is $2d+2$.
\end{exmp}

\begin{exmp} {\slsf (Odd degrees).} \rm
\label{exmp4}
Let $(a_1,b_1),...,(a_d,b_d)$ be such as in Example \ref{exmp3} and, in
addition, such that quadratic polynomials $(z-a_1)(z-b_1),...,(z-a_d)(z-b_d)$
are relatively prime with $(z-1)(z-2)$. Consider the rational map
\begin{equation}
f(z_1,z_2) = \Big( \frac{z_2}{z_1} , \frac{(z_1-1)(z_1-2)(3z_1-2)
\prod_{j=1}^d(z_1-a_j)(z_1-b_j)}
{(z_2-3z_1 +2)\prod_{j=1}^d[z_2-(a_j+b_j)z_1 + a_jb_j]}\Big).
\end{equation}
Polynomials $P(z_1,z_2) = (z_1-1)(z_1-2)(3z_1-2)\prod_{j=1}^d(z_1-a_j)(z_1-b_j)$
and $Q(z_1,z_2) = \prod_{j=1}^d[z_2-(a_j+b_j)z_1 + a_jb_j]$
are again relatively prime and
\[
P(z,z^2) - z^2Q(z,z^2) = (z-1)(z-2)(3z-2-z^2)\prod_{j=1}^d(z-a_j)(z-b_j)  =
\]
\[
= - (z-1)^2(z-2)^2\prod_{j=1}^d(z-a_j)(z-b_j),
\]
where $P(z) = (z-1)(z-2)(3z-2)\prod_{j=1}^d(z-a_j)(z-b_j)$ and $Q(z,z^2) =
(z-1)(z-2)\prod_{j=1}^d(z-a_j)(z-b_j)$.
Therefore the condition of Lemma \ref{lem-z-z2} is again satisfied. The
degree of $f$ is $2d+3$.
\end{exmp}

\newprg[prg2.3]{More examples}

One can produce in a similar way other examples. For example one may
take as $f_1(z_1,z_2) = \frac{z_2}{z_1^k}$ and then look for
$f_2(z_1,z_2) = \frac{P(z_1,z_2)}{Q(z_1,z_2)}$ with the condition
that $P$ and $Q$ are relatively prime and every non-zero root of
$P(z,z^{k+1}) - z^{k+1}Q(z,z^{k+1})$ should be a root of both
$P(z,z^{k+1})$ and $Q(z,z^{k+1})$. With such $f_2$ the map
$f=(f_1,f_2)$ will not have fixed points. Its degree will be at
least $k \text{ plus degree of } f_2 \text{ in } z_2$.

\begin{exmp} \rm
Consider the following map
\[
f : (z_1,z_2) \to \Big(\frac{z_2}{z_1^k}, \frac{z_1^{k+1} -1}{z_2-1}\Big).
\]
Then
\[
P(z,z^{k+1}) - z^{k+1}Q(z,z^{k+1}) = -(z^{k+1}-1)^2
\]
and therefore $f$ has no fixed points. Its degree is $k+1$.

\end{exmp}

\begin{exmp} \rm
One can start with $f_1(z_1,z_2)= \frac{z_2^{k}}{z_1^{k-1}}$. Then
the condition for a coprime $P, Q$ to define a map
$f=(f_1,\frac{P}{Q})$ without fixed points is this: for every
$k$-root of $1$ (denote it by $\zeta_l$, $l=1,...,k$), every
non-zero $z$, which satisfies $P(z,\zeta_l z) - \zeta_l zQ(z,\zeta_l
z)=0$, should also be the root of both $P(z,\zeta_l z)$ and
$Q(z,\zeta_l z)$. The following map will do the job:

\[
f: (z_1,z_2) \to \Big(\frac{z_2^{k}}{z_1^{k-1}}, \frac{z_1^k - 1}{z_2^k - 1}\Big).
\]
\end{exmp}

\newsect[sect3]{Transfer to the complex projective plane}

Now let us explain how to translate our examples from
$\pp^1\times\pp^1$ to $\pp^2$. This task is not completely obvious
because of rationality (\ie non-holomorphicity) of this transfer.
Fix some point $p=(p_1,p_2)\in\pp^1\times\pp^1$ and denote as
$l_1\deff \{p_1\}\times\pp^1$ and $l_2\deff \pp^1\times \{p_2\}$ the
``vertical'' and ``horizontal''  lines passing through $p$. Blow up
$X\deff\pp^1\times\pp^1$ at $p$ and then blow down the strict
transforms of $l_1$ and $l_2$. The obtained surface we denote as
$\hat X_p$. This $\hat X_p$ is biholomorphic to $\pp^2$. If $f:X\to
X$ is a rational self-map of $X$ then it naturally lifts to a
rational self-map $\hat f_p:\hat X_p\to\hat X_p$ of $\hat X_p$.

\begin{prop}
\label{transfert} Let $f:\pp^1\times\pp^1\to\pp^1\times\pp^1$ be a
rational map without fixed points. Take a regular point point
$p=(p_1,p_2)$ of $f$ such that:

\smallskip
\sli all $\{q_1,...,q_d\}\deff f^{-1}(p)$ and $s\deff f(p)$ are also
regular points of $f$;

\slii $l_1$ and $l_2$  do not intersect neither $\{s,q_1,...,q_d\}$
no $I(f)$ and are not contracted by $f$.

\noindent Then $\hat f_p$ has no fixed points to.
\end{prop}
\proof Denote by $E_p$ the exceptional divisor of the blow-up, by
$X_p$ the obtained surface and by $\pi_0 : X_p\to X$ the
corresponding blow-down map which sends $E_p$ to $\{p\}$. Lift $f$
to $X_p$ and denote by $f_p:X_p\to X_p$ the lifted map,  \ie
$f_p\deff \pi_0^{-1}\circ f\circ \pi_0$.

\medskip\noindent{\slsf Claim 1. {\it $f_p$ has no fixed points.}}
By $q_j$ denote the preimages of $q_j$ for $j=1,...,d$ and by $s$
that of $s$. Likewise denote by $I(f)$ the indeterminacy set of $f$
as well as its preimage under $\pi_0$. Then $I(f_p) = \{I(f);
q_1,..., q_d\}$. Really, $\pi_0$ is biholomorphic near every point
from $I(f)$ and therefore it remains an indeterminacy point also for
$f_p$. As for, say $q_1$, the map $f$ sends it to $p$ and
$\pi_0^{-1}$ blows it up. So $q_1$ becomes to be an indeterminacy
point of $f_p$. For any other point $r\in X_p$ both $\pi_0$ is
regular at $r$ and $f$ is regular at $\pi_0(r)$ and
$f(\pi_0(r))\not=p$. Therefore $\pi_0^{-1}$ is regular at
$f(\pi_0(r))$.

\smallskip Suppose $r\in X_p$ is a fixed point for $f_p$. Then
$r\not\in \{I(f); q_1,..., q_d\}$. If in addition $r\not\in E_p$
then both $r$ and $f_p(r)$ belong to the domain where $\pi_0$ is
biholomorphic. Therefore $f_p(r)=r$ would imply $f(\pi_0(r)) =
\pi_0(r)$ and this is not the case. The only case left is $r\in
E_p$. But hen $\pi_0$ sends $r$ to $p$ and $f$ further to $s\not=p$.
Finally $\pi^{-1}(s)\not\in E_p$ and we are done.

\medskip  Let $\pi_1:X_p\to X_1$ be the blow-down of $l_1$ and let
$f_1\deff \pi_1\circ f_p\circ \pi_1^{-1}$ be the pulled down map.
Set $s_1=\pi_1(l_1)$.

\medskip\noindent{\slsf Claim 2.  {\it $f_1$ has no fixed points
as well.}} Start form $s_1$, which is the only ``new'' point in
$X_1$. It is  an indeterminacy point for $f_1$. Really,
$\pi_1^{-1}(s_1)=l_1$ and $f_p(l_1)$ is not a point because $f(l_1)$
is not a point by the assumption. At the same time
$f_p(l_1)\not\subset E_p$, because again $f(l_1) = \pi_0(f(l_1))$
should not be a point. Therefore $s_1$ is not a fixed point of
$f_1$.

Take any other point $r\in X_1$. $\pi_1$ is biholomorphic near $r$
and therefore would $r$ be a fixed point for $f_1$ the point
$\pi_1^{-1}(r)$ would be a fixed point for $f_p$. Which cannot
happen according to the Case 1.

\medskip Let $\pi_2:X_1\to X_2$ be the blow-down of $l_2$ and let
$f_2\deff \pi_2\circ f_1\circ \pi_2^{-1}$ be the pulled down map.
Set $s_2=\pi_2(l_2)$.

\medskip\noindent{\slsf Claim 3. {\it $f_2$ has no fixed points to.}}
The proof is the same as for Claim 2. But $X_2\equiv \cc\pp^2$ and
we are done.

\smallskip\qed

\medskip The relevance of this Proposition to our task is clear: for
a given map $f$ a generic choice of $p$ will satisfy conditions
(\sli, (\slii of Proposition \ref{transfert} and therefore the
lifted map $\hat  f_p$ will be also without fixed points if such was
$f$.

\newsect[sect4]{Closure of the set of fixed point free rational
maps}

In Theorem \ref{fixed} we proved that the set of fixed point free
rational self-maps of $\pp^2$ is not empty. Now we are going to
prove that it is also closed modulo degenerate maps.

\newprg[prgCLOS.topol]{Topology on the space of rational maps}

Possible notions of convergence on the space of meromorphic mappings
$\calm (X,Y)$ between complex manifolds (or, spaces) are discussed
in \cite{I}. For us in this paper an appropriate one is the
following. Let $\{ f_n \}$ be a sequence of meromorphic maps from
complex manifold into $X$ a complex manifold $Y$.

\begin{defi}
\label{conv} We say that $f_n$  {\slsf strongly converges} on
compacts in $X$ to a meromorphic map $f : X \to Y$ if for any
compact $K\comp X$
\begin{equation}
\calh- \lim_{n\rightarrow \infty}\Gamma_{f_n}\cap (K\times Y) =
\Gamma_f\cap (K\times Y).
\end{equation}
\end{defi}
\noindent Here by $\calh-\lim $ we denote the limit in the Hausdorff
metric, supposing that both $X$ and $Y$ are equipped with some
Hermitian metrics. Remark that this notion of convergency doesn't
depend on a particular choice of metrics. We say  that $f_n$  {\slsf
converge weakly}, if their graphs converge i Hausdorff metric, see
\cite{I} for more details about relation between strongly and weakly
converging sequences of meromorphic mappings.

\smallskip The definition of strong convergence agrees well with the
usual notion of convergence of holomorphic mappings. Namely, in
\cite{I} we proved the following:

\begin{thm} {\slsf (Rouch\'e principle).} Let a sequence of meromorphic
mappings $\{ f_n\} $ between normal complex spaces $X$ and $Y$
converge strongly on compacts in $X$ to a meromorphic map $f$. Then:

\smallskip
\sli  If $f$ is holomorphic then for any relatively compact open
subset $X_1\comp X$ all restrictions $f_n\mid_{X_1}$ are holomorphic
for $n$ big enough, and $f_n \longrightarrow f$ on compacts in $X$
in the usual sense (of holomorphic mappings).

\slii If $f_n$ are holomorphic  then $f$ is also holomorphic and
$f_n \longrightarrow f$ on compacts in $X$.
\end{thm}

\newprg[prgCLOS.proof]{Proof of Theorem \ref{closed}}

In the case of compact $X$ the notion of strong convergence (and
corresponding topology on the space $\calm (X)$) is also well
adapted for understanding the structure of the space of meromorphic
self-maps of $X$. We shall not discuss here the general case any
more, but just turn to our case $X=\pp^2$ (or, $X=\pp^1\times
\pp^1$). The space $\calm (\pp^2)$ naturally splits into a double
sequence $\{\calm_{d,s}\}_{d=1,s=0}^{\infty}$ of subspaces, indexed
by the degree $d$ and the ``skew-degree'' $s$. Each $\calm_{d,s}$ is
infinite dimensional. If we equip $\pp^2$  with the Fubini-Studi
metric form $\omega$ of total volume one then
\begin{equation}
\eqqno(integrals)
d = \int\limits_{\pp^2}(f^*\omega)^2 \qquad \text{ and } \qquad s =
\int\limits_{\pp^2}\omega\wedge f^*\omega .
\end{equation}
The  volume of the graph of $f$ will be $\vol (\Gamma_f) = d+1+2s$.
Now if $\{f_n\}\subset\calm_{d,s}$, then by Bishop compactness
theorem we know that some subsequence $\{\Gamma_{f_{n_k}}\}$
converge to an analytic set $\Gamma$. This $\Gamma$ naturally
decomposes as the union of compact components: $\Gamma =
\Gamma_f\cup \Gamma_{\delta}\cup \Gamma_s$, where $\Gamma_f$ is a
graph of some meromorphic $f$, $\Gamma_{\delta}$ projects to points
and $\Gamma_s$ to curves under the natural projection
$\pi_1:\pp^2\times\pp^2\to \pp^2$ onto the source. Now $\vol
(\Gamma_{\delta})$ contributes to the first integral in
\eqqref(integrals), \ie to the degree and $\vol (\Gamma_{s})$ to the
second, \ie to the skew-degree. Therefore if the Hausdorff limit of
$\Gamma_{f_{n_k}}$ is different from $\Gamma_f$ our subsequence
diverge from $\calm_{d,s}$. Summing up we conclude that each
$\calm_{d,s}$ is closed in the strong topology we introduced (and is
relatively compact in the weak one).

\medskip Now let us prove the Theorem \ref{closed} from the Introduction.
Denote by $Fix (\pp^2)$ the set of rational self-mappings of $\pp^2$
admitting at least one fixed point. Let $f_n$ have no fixed points
and converge strongly to $f$. First of all we remark that Theorem 2
from \cite{I} implies that a strongly converging sequence has
bounded volume. From this fact and from what was said above we see
that $f_n$ for $n>>1$ and $f$ belong to the same component
$\calm_{d,s}$. If $f\in Fix(\pp^2)\setminus DFix(\pp^2)$ then
$\Gamma_f$ intersects the diagonal $D$ by an isolated point ${\bf
p}=(p,p)$. But then $f_n$, converging to $f$, should be holomorphic
in a neighborhood of $p$ and converge to $f$ as holomorphic
mappings by Rouch\'e's Principle. Therefore they should have their 
graphs intersect $D$ by the standard version of Rouch\'e theorem for 
vector functions. Contradiction.

\medskip Let us end up with an explicit example when a sequence without
fixed points converge to a map with a curve of fixed points.

\begin{exmp} \rm Take $\theta_n$ to be irrational modulo $2\pi$ and
converging to zero. Take the following self-maps of $\pp^2$:
$f_n:[z_0:z_1:z_2] \to \Big[z_0z_1:z_0z_2:e^{i\theta_n}z_1z_2\Big]$.
Then  it is straightforward to check that $f_n$ have no fixed
points. But the limit map $f_n:[z_0:z_1:z_2] \to \Big[z_0z_1:z_0z_2:
z_1z_2\Big]$ has the  curve $\{ z_2z_0=z_1^2\}$ of fixed points.
\end{exmp}

\ifx\undefined\bysame
\newcommand{\bysame}{\leavevmode\hbox to3em{\hrulefill}\,}
\fi

\def\entry#1#2#3#4\par{\bibitem[#1]{#1}
{\textsc{#2 }}{\sl{#3} }#4\par\vskip2pt}

\end{document}